\documentstyle[11pt]{article}
\begin{document}
\setlength{\baselineskip}{17pt}
\begin{center}
{\LARGE {\bf Stability Conditions and Liapunov Functions \\
for Quasipolynomial Systems }}
\end{center}

\mbox{}

\begin{center}
{\large {\sc Benito Hern\'{a}ndez--Bermejo$^1$ }} \\
Departamento de F\'{\i}sica Matem\'{a}tica y Fluidos, Facultad de Ciencias \\
Universidad Nacional de Educaci\'{o}n a Distancia \\
Senda del Rey S/N, 28040 Madrid, Spain. \\
\end{center}

\mbox{}

\mbox{}

\mbox{}

\noindent{\bf Abstract---}
The stability of equilibrium points of quasipolynomial systems of ODEs is 
considered. The criteria and Liapunov functions found generalize those 
traditionally known for Lotka-Volterra equations, that now appear as a 
particular case. 

\mbox{}

\mbox{}

\mbox{}


\mbox{}

\mbox{}

\mbox{}

\mbox{}

\mbox{}

\mbox{}

\mbox{}

\mbox{}

\mbox{}

\mbox{}

\mbox{}

\mbox{}

\mbox{}

\noindent $^1$ E-mail: bhernand@apphys.uned.es \newline
Telephone: (+ 34 91) 398 72 19. \hspace{3mm} Fax: (+ 34 91) 398 66 97.
\pagebreak

Consider the well-known Lotka-Volterra (LV) system 
\begin{equation}
   \label{lv0}
   \dot{x}_i = x_i \left( \lambda _i + \sum _{j=1}^n A_{ij} x_j \right) 
   , \;\:\;\: i = 1, \ldots ,n
\end{equation}
which is assumed to have a unique positive equilibrium point: 
$x^* = (x_1^*, \ldots ,x_n^*)^T \in$ int$( I\!\!R_+^n)$. A question of 
fundamental importance in the analysis of these equations concerns the 
stability of the steady state (see [\ref{gh1}, pp. 3-5] for the main 
definitions of stability, which are the ones adopted in this work). In this 
sense, many of the most relevant results about stability of LV systems are 
those based on the Liapunov function [\ref{goh1}--\ref{tak1}]: 
\begin{equation}
    \label{hlv0}
    V(x) = \sum _{i=1}^{n} c_i \left( x_i-x_i^* -
      x_i^* \ln \frac{x_i}{x_i^*} \right) \;\; , \;\:\; c_i>0 \; \forall \; i 
\end{equation}
The time derivative of (\ref{hlv0}) along the trajectories of (\ref{lv0}) is
\begin{equation}
    \label{hp0}
    \dot{V}(x) = \frac{1}{2}(x-x^*)^T(C \cdot A+A^T \cdot C)(x - x^*)
\end{equation}
where $C=$ diag$(c_1, \ldots ,c_n)$. Thus, it can be stated that if there 
exists a positive definite diagonal matrix $C$ such that 
$C \cdot A + A^T \cdot C$ is negative semi-definite (i.e. $A \in \bar{S}_w$) 
then $x^*$ is stable. Moreover if, instead, $C \cdot A + A^T \cdot C$ is 
negative definite (that is, $A \in S_w$) then $x^*$ is globally 
asymptotically stable in int$(I \!\! R^n_+)$. 

A natural and extensive generalization of LV systems is provided by 
quasipolynomial (QP) systems of ODEs (see [\ref{bvb1}] and references 
therein):
\begin{equation}
\label{qps}
    \dot{x}_i = x_i \left( \lambda_i + \sum_{j=1}^m A_{ij} 
    \prod_{k=1}^n x_k^{B_{jk}} \right) \; , \;\:\; i = 1, \ldots , n \; , 
    \;\:\; m \geq n \; , \;\:\; \mbox{Rank}(B) = n
\end{equation}
Conditions $m \geq n$ and Rank($B$) $= n$ can be assumed without loss of 
generality [\ref{bvb1}] because systems not verifying them are reducible to 
the form (\ref{qps}) in int$(I \!\! R^n_+)$. Therefore, (\ref{qps}) is a 
completely generic starting point for what is to follow. 

The aim of this letter is to demonstrate that the aforementioned properties 
regarding the stability of LV equations (\ref{lv0}) can be generalized to a 
great extent for QP systems (\ref{qps}). This is interesting from an applied 
perspective, because most systems of ODEs arising in practice are QP or can 
be algorithmically recasted in such form [\ref{bvb1},\ref{bv1}].

For convenience let us first define the $m \times m$ matrix $Q = B \cdot A$, 
as well as the quasimonomial functions:
\begin{equation}
\label{qms}
    \varphi _i(x)= \prod_{j=1}^n x_j^{B_{ij}} \; , \;\:\; i = 1, \ldots , m 
\end{equation}

The main result is the following one: 

\mbox{}

\noindent{\sc Theorem 1.} {\sl For every fixed point $x^* \in$ 
int$(I \!\! R^n_+)$ of a QP system (\ref{qps}):

\noindent a) If $Q \in \bar{S}_w$ then $x^*$ is stable.

\noindent b) If $Q \in S_w$ then $x^*$ is globally asymptotically stable in 
int$(I \!\! R^n_+)$.

\noindent c) Let $C=$ diag$(c_1, \ldots ,c_m)$, with $c_i>0 \; \forall \; i$, 
be a matrix such that $C \cdot Q + Q^T \cdot C$ is negative semi-definite in 
case (a) or negative definite in case (b). Then, in either (a) or (b), 
\begin{equation}
\label{liap}
   W(x) = \sum_{i=1}^m c_i \left\{ 
   \prod_{j=1}^n x_j^{B_{ij}} - 
   \prod_{j=1}^n (x^*_j)^{B_{ij}} - \left( \prod_{j=1}^n (x^*_j)^{B_{ij}} 
   \right) \ln \left[ \prod_{j=1}^n \left( \frac{x_j}{x_j^*} \right)^{B_{ij}} 
   \right] \right\}
\end{equation}
is a Liapunov function for $x^*$ in int($I \!\! R^n_+$).

\noindent d) The time variation of the Liapunov function (\ref{liap}) along 
the trajectories of (\ref{qps}) is:
\begin{equation}
\label{tvar}
   \dot{W} = \frac{1}{2} (\varphi(x) - \varphi(x^*))^T 
   (C \cdot Q + Q^T \cdot C)(\varphi(x) - \varphi(x^*))
\end{equation}
where $\varphi(x) = (\varphi_1(x) , \ldots , \varphi_m(x))^T$ and the 
$\varphi_i(x)$ are defined in (\ref{qms}).
}

\mbox{}

\noindent{\sc Proof.}
Notice that $W(x^*)=0$. In addition, we shall first demonstrate that 
function $W(x)$ in (\ref{liap}) verifies $W(x)>0$ for all $x \neq x^*$, 
$x \in$ int($I \!\! R^n_+$). For this, consider the $m$-dimensional function:
\begin{equation}
    \label{hlv1}
    V(y_1, \ldots ,y_m) = \sum _{i=1}^{m} c_i \left( y_i-y_i^* -y_i^* \ln 
    \frac{y_i}{y_i^*} \right) \;\; , \;\:\; c_i>0 \; \forall \; i
    \;\; , \;\:\; y_i^*>0 \; \forall \; i
\end{equation}
Obviously, $V(y^*)=0$ and $V(y)>0$ for every $y \neq y^*$, $y \in$ 
int($I \!\! R^m_+$). Let us now perform the following change of variables:
\begin{equation}
\label{qmt}
   y_i = \prod_{j=1}^m z_j^{\tilde{B}_{ij}} \; , \;\:\; i = 1, \ldots , m
\end{equation}
In (\ref{qmt}) we define $\tilde{B} = \left(B \mid B^*_{m \times (m-n)} \right)$, 
where $B^*$ is an $m \times (m-n)$ submatrix of arbitrary entries to be 
chosen in such a way that $\tilde{B}$ is invertible (note that this is always 
possible since Rank($B$) $=n$). Consequently, transformation (\ref{qmt}) is 
bijective in int($I \!\! R^m_+$). Let $z^* =(z_1^*, \ldots , z_m^*)^T \in$ 
int($I \!\! R^m_+$) be the unique point whose image is $y^*$ after 
transformation (\ref{qmt}). Then, function $V(y)$ in (\ref{hlv1}) is mapped 
to the $m$-dimensional function $\tilde{W}(z) = V(y(z))$, defined in 
int($I \!\! R^m_+$): 
\begin{equation}
\label{tliap}
   \tilde{W}(z) = \sum_{i=1}^m c_i \left\{ 
   \prod_{j=1}^m z_j^{\tilde{B}_{ij}} - 
   \prod_{j=1}^m (z^*_j)^{\tilde{B}_{ij}} - 
   \left( \prod_{j=1}^m (z^*_j)^{\tilde{B}_{ij}} \right) 
   \ln \left[ \prod_{j=1}^m \left( \frac{z_j}{z_j^*} 
   \right)^{\tilde{B}_{ij}} \right] \right\}
\end{equation}
Thus we have that $\tilde{W}(z^*)=0$ and $\tilde{W}(z) > 0$ for every 
$z \neq z^*$, $z \in$ int($I \!\! R^m_+$). 

If $m=n$, then $\tilde{B}$ reduces to $B$ and if we replace 
in (\ref{tliap}) $z$ and $z^*$ by $x$ and $x^*$, respectively, it is 
demonstrated that $W(x)$ in (\ref{liap}) verifies $W(x)>0$ for all 
$x \neq x^*$, $x \in$ int($I \!\! R^n_+$). 

To prove the same in the complementary $m>n$ case, let us choose a point 
$z^*$ of the form 
$z^* =(x_1^*, \ldots ,x_n^*, 1, \stackrel{(m-n)}{\ldots} ,1)^T$. We know that 
$\tilde{W}(z)>0$ for every $z \neq z^*$, $z \in$ int($I \!\! R^m_+$). Then 
this will be the case, in particular, for the points of the hyperplane 
$\{ z_i=1, i = n+1, \ldots, m \}$ belonging to int($I \!\! R^m_+$). Let us 
parametrize those points as $z=(x_1, \ldots ,x_n, 1, \stackrel{(m-n)}{\ldots} 
,1)^T$, with $x_i>0$ for all $i=1, \ldots ,n$. Thus if we substitute $z^*$ and 
the parametrized form of $z$ in $\tilde{W}$, it is immediately obtained that 
$\tilde{W}(x_1, \ldots ,x_n, 1, \stackrel{(m-n)}{\ldots},1) = W(x_1, \ldots , 
x_n)$ in (\ref{liap}) because $\tilde{B}_{ij}=B_{ij}$ for $j=1, \ldots, n$. 
This proves that $W(x)>0$ for all $x \neq x^*$, $x \in$ int($I \!\! R^n_+$), 
in the $m>n$ case. 

Finally, let us look at the time evolution of $W(x)$. For this, note that 
we can write the QP system (\ref{qps}) as:
\begin{equation}
\label{qps2}
    \dot{x}_i = x_i \sum_{j=1}^m A_{ij} ( \varphi_j(x) - \varphi_j(x^*) )
    \; , \;\:\; i = 1, \ldots , n 
\end{equation}
where $\varphi_i(x)$ is defined in (\ref{qms}). Note also that 
\begin{equation}
\label{parf}
   \frac{\partial \varphi _i (x)}{\partial x_j} = \frac{B_{ij} 
   \varphi _i (x)}{x_j}
\end{equation}
We then have:
\begin{equation}
   \dot{W} = 
   \sum_{i=1}^m c_i \dot{\varphi}_i(x) \left( 1 - 
   \frac{\varphi_i(x^*)}{\varphi_i(x)} \right) = 
   \sum_{i=1}^m c_i \left( 1 - \frac{\varphi_i(x^*)}{\varphi_i(x)} \right) 
   \sum_{k=1}^n \frac{\partial \varphi _i (x)}{\partial x_k} \dot{x}_k
\end{equation}
Taking (\ref{qps2}) and (\ref{parf}) into account, (\ref{tvar}) is found 
after some simple algebra. 

Consequently [\ref{gh1}, pp. 3-5], we have that if $Q \in \bar{S}_w$ 
($Q \in S_w$) then $W(x)$ is a Liapunov function for the system in 
int($I \!\! R^n_+$) and $x^*$ is stable (globally asymptotically stable in 
int($I \!\! R^n_+$)). 

This completes the proof of the theorem.

\mbox{}

\noindent{\sc Remark 1.} Notice that the aforementioned results for LV 
systems appear now as a particular case of Theorem 1 when $m=n$, $B$ is the 
identity matrix and Rank($A$) $=n$. 

\mbox{}

\noindent{\sc Remark 2.} The previously known criteria for the belonging of 
a matrix to $S_w$ or $\bar{S}_w$, widely investigated in the context of the 
stability of LV equations (see [\ref{tak1}] and references therein) can now 
be extended straightforwardly to QP systems.

\mbox{}

\noindent{\sc Corollary 1.} {\sl Under the same hypotheses of Theorem 1, if 
$\dot{W}(x)=0$ for all $x \in$ int$(I \!\! R^n_+)$, then the trajectories of 
system (\ref{qps}) lie on the surfaces:}
\begin{equation}
\label{ham}
   \sum_{i=1}^m c_i \left\{ 
   \prod_{j=1}^n x_j^{B_{ij}} - 
   \prod_{j=1}^n (x^*_j)^{B_{ij}} - \left( \prod_{j=1}^n (x^*_j)^{B_{ij}} 
   \right) \ln \left[ \prod_{j=1}^n \left( \frac{x_j}{x_j^*} \right)^{B_{ij}} 
   \right] \right\} = \mbox{constant}
\end{equation}

\mbox{}

\noindent{\sc Remark 3.} Corollary 1 takes place, for instance, in QP systems 
described in terms of a Poisson structure, which are closely related 
[\ref{pst}] to conservative LV systems (see [\ref{log1}] for a review of LV 
conservativeness). In such cases the Liapunov function (\ref{ham}) is a first 
integral playing the role of Hamiltonian.

\mbox{}

\mbox{}

\begin{center}
{\bf REFERENCES}
\end{center}
\begin{enumerate}
\item J. Guckenheimer and P. Holmes, {\em Nonlinear Oscillations, Dynamical 
      Systems, and Bifurcations of Vector Fields,\/} Springer-Verlag, New 
      York (1983).\label{gh1}
\addtolength{\itemsep}{-2mm}
\item B. S. Goh, {\em Management and Analysis of Biological Populations,\/} 
      Elsevier, Oxford-Amsterdam-New York (1980).\label{goh1}
\item J. Hofbauer and K. Sigmund, {\em The Theory of Evolution and Dynamical 
      Systems,\/} Cambridge University Press, Cambridge (1988).\label{hs1}
\item D. O. Logofet, {\em Matrices and Graphs. Stability Problems in 
      Mathematical Ecology,\/} CRC Press, Boca Raton, Florida 
      (1993).\label{log1}
\item Y. Takeuchi, {\em Global Dynamical Properties of Lotka-Volterra 
      Systems,\/} World Scientific, Singapore (1996).\label{tak1}
\item B. Hern\'{a}ndez--Bermejo, V. Fair\'{e}n and L. Brenig, Algebraic 
      recasting of nonlinear systems of ODEs into universal formats, {\em J. 
      Phys. A, Math. Gen.\/} {\bf 31}, 2415-2430 (1998).\label{bvb1}
\item B. Hern\'{a}ndez--Bermejo and V. Fair\'{e}n, Lotka-Volterra 
      representation of general nonlinear systems, {\em Math. Biosci.\/} 
      {\bf 140}, 1-32 (1997).\label{bv1}
\item B. Hern\'{a}ndez--Bermejo and V. Fair\'{e}n, Hamiltonian structure and 
      Darboux theorem for families of generalized Lotka-Volterra systems, 
      {\em J. Math. Phys.\/} {\bf 39}, 6162-6174 (1998).\label{pst}
\end{enumerate}
\end{document}